\documentstyle[a4,12pt]{article}
\def\beq{\begin{equation}}
\def\eeq{\end{equation}}
\def\bea{\begin{eqnarray}}
\def\eea{\end{eqnarray}}
\def\nn{\nonumber}
\def\jor{$ {\cal U}_h(sl(2)) \ $}
\def\jorsu{$ {\cal U}_h(su(1,1)) \ $}
\def\X2{{hX \over 2}}
\def\Zp2{{hZ_+ \over 2}}
\def\R2{{hR \over 2}}
\def\Tp2{{hT_+ \over 2}}

\def\ket#1{\left| #1\right\rangle}
\begin{document}
\thispagestyle{empty}
\begin{flushright}
OWUAM-020 \\
January 20, 1997
\end{flushright}

\begin{center}

\vfill
{\large Irreducible Decomposition for Tensor Product Representations of 
Jordanian Quantum Algebras}
\vspace{1.5cm}

N. Aizawa
\vspace{0.5cm}

{\em Department of Applied Mathematics}

{\em Osaka Women's University, Sakai, Osaka 590, Japan}

\end{center}

\vfill
\begin{abstract}
  Tensor products of irreducible representations of the Jordanian 
quantum algebras \jor and \jorsu are considered. For both the highest 
weight finite dimensional representations of \jor and lowest 
weight infinite dimensional ones of \jorsu, 
it is shown that tensor product representations are reducible 
and that the decomposition rules to irreducible representations are 
exactly the same as those of corresponding Lie algebras. 
\end{abstract}
\newpage

\setcounter{equation}{0}
\section{Introduction}

  Recent works on quantum matrices in two dimensions \cite{dmmz,eow} 
introduced a new deformation of the Lie algebra $sl(2)$ 
called $h$-deformation or Jordanian deformation \jor \cite{ohn}. 
Some algebraic aspects of \jor have been investigated and 
it has been shown that \jor is a quasitriangular Hopf algebra 
\cite{sak,bh} and that \jor can be constructed from the Drinfelf-Jimbo 
deformation by a contraction \cite{aks}. Furthermore two kinds of mappings 
from $sl(2)$ to \jor have been obtained \cite{kob,acc1}. 

  On the other hand, representation theories of \jor have not been 
welldeveloped yet. What we know so far is that the finite dimensional 
irreducible representations of \jor are classified exactly the same 
way as those of $sl(2)$. To show this, the standard singular vector 
construction method is used in \cite{dov,acc2}, while 
the authors of \cite{acc1} and \cite{bfhn} use the nonlinear invertible 
map from $sl(2)$ to \jor and boson realizations, respectively. 
In \cite{bfhn}, it is shown that decomposition rules of tensor 
product representations are the same as $sl(2)$ for some low 
dimensional representations.

  In this paper, we consider irreducible decomposition for 
tensor product representations of Jordanian quantum algebras. 
Representations discussed in this paper are highest weight 
finite dimensional ones for \jor and lowest weight infinite 
dimensional ones for \jorsu. The Jordanian quantum algebra 
\jorsu is introduced as an algebra being isomorphic to \jor. 
It is shown that the decomposition rules for both cases are 
the same as the classical cases. Some examples are shown for 
\jor in order to discuss explicit expressions of 
Clebsch-Gordan coefficients. 
This work is motivated by the fact that 
welldeveloped representation theories are necessary when we consider 
physical applications of algebraic objects.

\setcounter{equation}{0}
\section{\jor and its representations}

  The Jordanian quantum algebra \jor is an associative algebra with 1 
generated by $ X, Y $ and $ H$. Their commutation relations are 
given by \cite{ohn}
\bea
& & [H, \; X] = 2{\sinh hX \over h}, \qquad 
    [H, \; Y] = -Y (\cosh hX) - (\cosh hX) Y, \nn \\
& & [X, \; Y] = H,                            \label{com1}
\eea
where $h$ is the deformation parameter. The Casimir element is 
\beq
C = {1 \over 2h} \left\{ Y(\sinh hX) + (\sinh hX)Y \right\} + {1 \over 4} H^2 
  + {1 \over 4} (\sinh hX)^2. 
                                               \label{cas}
\eeq
In the limit of $ h \longrightarrow 0 $, \jor reduces to $sl(2)$. 
The Hopf algebra structure reads
\bea
 & & \Delta(X) = X \otimes 1 + 1 \otimes X, \nn \\
 & & \Delta(Y) = Y \otimes e^{hX} + e^{-hX} \otimes Y, \nn \\
 & & \Delta(H) = H \otimes e^{hX} + e^{-hX} \otimes H, \label{hopf} \\
 & & \epsilon(X) = \epsilon(Y) = \epsilon(H) = 0, \nn \\
 & & S(X) = -X, \quad S(Y) = -e^{hX}Ye^{-hX}, \quad 
     S(H) = -e^{hX} H e^{-hX}.                \nn
\eea

  The finite dimensional highest weight representations can be 
easily obtained by making use of the invertible map from $sl(2)$ 
to \jor given in \cite{acc1}. Let us define the following elements 
according to \cite{acc1}
\bea
 & & Z_+ = {2 \over h} \tanh \X2, \nn \\
 & & Z_- = (\cosh \X2) Y (\cosh \X2), \label{rl}
\eea
then it is not difficult to verify directly that $ Z_{\pm}$ and $H$ satisfy 
the $sl(2)$ commutation relations 
\beq
  [H, \; Z_{\pm}] = \pm 2Z_{\pm}, \qquad 
  [Z_+,\; Z_-] = H,               \label{sl2}
\eeq
and the Casimir element reads
\beq
  C = Z_+ Z_- + {H \over 2} \left({H \over 2} -1 \right)_,  \label{cas2}
\eeq
by making use of the identities proved by the mathematical induction
\bea
 & & [H, \; X^n] = 2nX^{n-1}{\sinh hX \over h}, \nn \\
 & & [Y, \; X^n] = -nX^{n-1}H - n(n-1)X^{n-2}{\sinh hX \over h}, 
     \label{beki}
\eea
where $n$ is a natural number. The authors of \cite{acc1} regard 
$Z_{\pm}, H$ as elements of $sl(2)$, however it is more convenient to 
regard them as elements of \jor for our purpose. Namely, their coproducts 
are given in terms of $ \Delta(X), \Delta(Y)$ and $ \Delta(H)$. 

 From (\ref{sl2}) and (\ref{cas2}), it is obvious that we can 
take the following as the irreducible highest weight  
representations of \jor
\bea
 & & Z_+ \ket{j\; m} = \ket{j\; m+1}, \nn \\
 & & Z_- \ket{j\; m} = (j+m)(j-m+1) \ket{j\; m-1}, \label{sl2rep} \\
 & & H \ket{j\; m} = 2m \ket{j\; m}, \nn 
\eea
and the eigenvalues of the Casimir element is
\beq
 C \ket{j\; m} = j(j+1) \ket{j\; m},         \label{caseigen}
\eeq
where $ j$ is a halfinteger or an integer and 
$ m= -j, -j+1, \cdots j.$ 
We adopt the unfamiliar representation for physicists for the 
sake of simplicity of calculations. This choice of the 
representations is not essential. All the discussions in subsequent 
sections hold for the usual representations. The representation matrices 
for $ X, \ Y $ can be obtained by solving (\ref{rl}) with respect to 
$ X,\ Y $ \cite{acc1}.

\setcounter{equation}{0}
\section{Eigenvectors of $ \Delta(H) $}

  We consider irreducible decomposition of tensor product of 
two representations given by (\ref{sl2rep}) ; 
$ \ket{j_1\; m_1} \otimes \ket{j_2\; m_2}. $ 
The key of deriving a decomposition rule is to construct the 
eigenvectors of $ \Delta(H) $, since if we obtain such vectors, 
the decomposition rules can be derived by the same discussion 
as the case of $sl(2)$ as we shall see later.

  First, we rewrite $ \Delta(H)$ in terms of $ H$ and $Z_{\pm}$. 
From (\ref{rl})
\bea
& & e^{hX} = {1 + \Zp2 \over 1 - \Zp2} 
    = 1 + 2\sum_{n=1}^{\infty}\; \left(\Zp2 \right)^n_, \nn \\
& & e^{-hX} = {1 - \Zp2 \over 1 + \Zp2} 
    = 1 + 2\sum_{n=1}^{\infty}\; \left(-\Zp2 \right)^n_, 
                                                   \label{ehx}
\eea
we obtain
\beq
 \Delta(H) = H \otimes 1 + 1 \otimes H 
  + H \otimes 2\sum_{n=1}^{\infty}\; \left(\Zp2 \right)^n 
  + 2\sum_{n=1}^{\infty}\; \left(-\Zp2 \right)^n \otimes H.
                                                  \label{coproh}
\eeq

  We denote an eigenvector of $ \Delta(H) $ whose eigenvalue 
is $2(m_1+m_2)$ by $ \ket{(j_1m_1)\; (j_2m_2)}. $ From (\ref{coproh}), 
$ \ket{(j_1m_1)\; (j_2m_2)} $ may be written as
\beq
  \ket{(j_1m_1)\; (j_2m_2)} = \sum_{k=0}^{j_1-m_1} \sum_{l=0}^{j_2-m_2} \; 
  \alpha(m_1+k,m_2+l) \ket{j_1 \; m_1+k} \otimes \ket{j_2 \; m_2+l}. 
                                                  \label{ansatz}
\eeq
We take $ \alpha(m_1, m_2) = 1 $ so as to reduce to the correct limit of 
$ h \longrightarrow 0.$ 
Substituting (\ref{coproh}) and (\ref{ansatz}) into 
\beq
   \Delta(H) \ket{(j_1m_1)\; (j_2m_2)} = 2(m_1+m_2) 
   \ket{(j_1m_1)\; (j_2m_2)},                     \label{eigeneq}
\eeq
we obtain 
\bea
 & & \sum_{k=0}^{j_1-m_1} \sum_{l=0}^{j_2-m_2} \;
  \{ (k+l) \alpha(m_1+k, m_2+l) + 2(m_1+k) \sum_{n=1}^l \; 
  \left({h \over 2}\right)^n \alpha(m_1+k, m_2+l-n)     \nn \\
 & & + 2(m_2+l) \sum_{n=1}^k \; \left(-{h \over 2}\right)^n 
 \alpha(m_1+k-n,m_2+l) 
 \} 
 \ket{j_1 \; m_1+k} \otimes \ket{j_2 \; m_2+l} =0. \label{sub}
\eea
Therefore, $ \alpha(m_1+k, m_2+l)$ must satisfy the recurrence relation
\bea
 & & (k+l) \alpha(m_1+k, m_2+l) + 2(m_1+k) \sum_{n=1}^l \; 
  \left({h \over 2}\right)^n \alpha(m_1+k, m_2+l-n)     \nn \\
 & & + 2(m_2+l) \sum_{n=1}^k \; \left(-{h \over 2}\right)^n 
 \alpha(m_1+k-n,m_2+l) =0.                      \label{rec1}
\eea

  Next, we rewrite the recurrence relation (\ref{rec1}) into a 
simpler form. Multiplying (\ref{rec1}) by $ -h/2 $ and replacing $k$ with 
$k-1$ then subtracting from (\ref{rec1}), the obtained relation reads
\bea
 & & (k+l) \alpha(m_1+k, m_2+l) - {h \over 2} (2m_2+1-k+l) 
     \alpha(m_1+k-1, m_2+l) \nn \\
 & & + 2\sum_{n=1}^l \; \left({h \over 2}\right)^n
     \{ (m_1+k) \alpha(m_1+k, m_2+l-n) 
     + {h \over 2} (m_1+k-1) \alpha(m_1+k-1, m_2+l-n) \} \nn \\
 & &                                      = 0. \label{rec2}
\eea
Multiplying (\ref{rec2}) by $h/2$ and replacing $l$ with $l-1$, then 
subtracting from (\ref{rec2}), we obtain the simpler form of 
recurrence relation
\bea
 & & (k+l) \; \alpha(m_1+k, m_2+l) + {h \over 2} (2m_1+1+k-l) 
     \; \alpha(m_1+k, m_2+l-1)                                  \nn \\
 & & - {h \over 2} (2m_2+1-k+l) \; \alpha(m_1+k-1,m_2+l)  \label{rec3} \\
 & & +\left({h \over 2}\right)^2 (2m_1+2m_2-2+k+l) \; 
      \alpha(m_1+k-1, m_2+l-1) = 0. \nn
\eea

  The solutions of the recurrence relation (\ref{rec3}) are given by
\bea
 & & \alpha(m_1+k, m_2+l) \nn \\
 & & = (-1)^l \left({h \over 2}\right)^{k+l} 
 \sum_{p=0}\; \left(
 \begin{array}{c}
  2m_1+k-p \\ l-p 
 \end{array}\right) 
 \left(
 \begin{array}{c}
 2m_1+k-1 \\ p
 \end{array}\right) 
 \left(
 \begin{array}{c}
 2m_2 \\ k-p
 \end{array}\right)_, 
 \label{sol}
\eea
where the sum on $p$ runs as far as all the binomial coefficients are 
welldefined. For the  negative values of $m_i$, the binomial coefficients 
are rewritten by the formula
\beq
  \left(
  \begin{array}{c}
    m \\ l
  \end{array}
  \right) 
  = (-1)^l
  \left( 
  \begin{array}{c}
    |m|+l-1 \\ l
  \end{array} \right)_.                  \label{negative}
\eeq
Substituting (\ref{sol}) into (\ref{rec3}), it can be verified that 
(\ref{sol}) gives the solutions of the recurrence relation (\ref{rec3}). 
We briefly sketch the calculation, since it is somewhat complicated. 
Substituting (\ref{sol}) into (\ref{rec3}), then using the identities
\bea
 & & (2m_1+1+k-l) \left( 
 \begin{array}{c}
  2m_1+k-p \\ l-1-p
 \end{array} \right) 
 =
 (l-p) \left(
 \begin{array}{c}
  2m_1+k-p \\ l-p
 \end{array} \right)_, \nn \\
 & & (2m_2+1-k+p) \left(
 \begin{array}{c}
  2m_2 \\ k-1-p
 \end{array} \right) 
 = 
 (k-p) \left(
 \begin{array}{c}
  2m_2 \\ k-p
 \end{array} \right)_,  \nn 
\eea
the left hand side of the recurrence relation (\ref{rec3}) can 
be rewritten 
\bea
 & & (-1)^l \left({h \over 2}\right)^{k+l} \left\{
 k \sum_{p=0}\; 
 \left( \begin{array}{c} 2m_1+k-p \\ l-p \end{array} \right) 
 \left( \begin{array}{c} 2m_1+k-1 \\ p \end{array} \right) 
 \left( \begin{array}{c} 2m_2 \\ k-p \end{array} \right) 
 \right. \nn \\
 & & 
 - \sum_{p=0}\; (k-p) 
 \left( \begin{array}{c} 2m_1+k-1-p \\ l-p \end{array} \right) 
 \left( \begin{array}{c} 2m_1+k-2 \\ p \end{array} \right) 
 \left( \begin{array}{c} 2m_2 \\ k-p \end{array} \right)  \nn \\
 & & 
 - \sum_{p=0}\; (l-p) 
 \left( \begin{array}{c} 2m_1+k-1-p \\ l-p \end{array} \right) 
 \left( \begin{array}{c} 2m_1+k-2 \\ p \end{array} \right) 
 \left( \begin{array}{c} 2m_2 \\ k-1-p \end{array} \right) \nn \\
 & &
 + \sum_{p=0}\; p 
 \left( \begin{array}{c} 2m_1+k-p \\ l-p \end{array} \right) 
 \left( \begin{array}{c} 2m_1+k-1 \\ p \end{array} \right) 
 \left( \begin{array}{c} 2m_2 \\ k-p \end{array} \right) \nn \\
 & & 
 - \sum_{p=0} \; (k-p) 
 \left( \begin{array}{c} 2m_1+k-1-p \\ l-1-p \end{array} \right) 
 \left( \begin{array}{c} 2m_1+k-2 \\ p \end{array} \right) 
 \left( \begin{array}{c} 2m_2 \\ k-p \end{array} \right) \nn \\
 & & 
 \left. 
 - \sum_{p=0} \; (2m_1-3+2k+l-p) 
 \left( \begin{array}{c} 2m_1+k-1-p \\ l-1-p \end{array} \right) 
 \left( \begin{array}{c} 2m_1+k-2 \\ p \end{array} \right) 
 \left( \begin{array}{c} 2m_2 \\ k-1-p \end{array} \right)
 \right\}_. \nn
\eea
Redefining $p+1$ as $p$ in the third and the sixth summation, the fourth 
and the sixth summation can be combined. The second and the fifth 
summation can also be combined by using the identity
\[
 \left( \begin{array}{c} n \\ l-1 \end{array} \right) 
 + 
 \left( \begin{array}{c} n \\ l \end{array} \right) 
 = 
 \left( \begin{array}{c} n+1 \\ l \end{array} \right)_.
\]
At this stage, the left hand side of (\ref{rec3}) reads
\bea
 & & (-1)^l \left({h \over 2}\right)^{k+l} \left\{
 k \sum_{p=0}\; 
 \left( \begin{array}{c} 2m_1+k-p \\ l-p \end{array} \right) 
 \left( \begin{array}{c} 2m_1+k-1 \\ p \end{array} \right) 
 \left( \begin{array}{c} 2m_2 \\ k-p \end{array} \right) 
 \right. \nn \\
 & & 
 - \sum_{p=0}\; (k-p) 
 \left( \begin{array}{c} 2m_1+k-p \\ l-p \end{array} \right) 
 \left( \begin{array}{c} 2m_1+k-2 \\ p \end{array} \right) 
 \left( \begin{array}{c} 2m_2 \\ k-p \end{array} \right)  \nn \\
 & & 
 - \sum_{p=1}\; (l-p+1) 
 \left( \begin{array}{c} 2m_1+k-p \\ l+1-p \end{array} \right) 
 \left( \begin{array}{c} 2m_1+k-2 \\ p-1 \end{array} \right) 
 \left( \begin{array}{c} 2m_2 \\ k-p \end{array} \right)  \nn \\
 & & \left.
 - \sum_{p=1}\; (k+l-1-p) 
 \left( \begin{array}{c} 2m_1+k-p \\ l-p \end{array} \right) 
 \left( \begin{array}{c} 2m_1+k-2 \\ p-1 \end{array} \right) 
 \left( \begin{array}{c} 2m_2 \\ k-p \end{array} \right)
 \right\}_.  \nn
\eea
It is now easy to see that this always vanishes, noting that the last 
two summation are combined to give
\[
   \sum_{p=1}\; (2m_1+2k-1-p) 
 \left( \begin{array}{c} 2m_1+k-p \\ l-p \end{array} \right) 
 \left( \begin{array}{c} 2m_1+k-2 \\ p-1 \end{array} \right) 
 \left( \begin{array}{c} 2m_2 \\ k-p \end{array} \right)_.  \nn \\
\]

  We therefore have shown that, for given vectors 
$ \ket{j_1m_1} $ and $ \ket{j_2m_2} $, a unique eigenvector of 
$ \Delta(H) $ with eigenvalu $ 2 (m_1 + m_2) $ 
can be constructed. The vector is given by 
(\ref{ansatz}) with $ \alpha(m_1+k, m_2+l) $ given by 
(\ref{sol}).

\setcounter{equation}{0}
\section{Decomposition rule for \jor}

  It has been shown in the previous section that we can construct 
a unique vector $ \ket{(j_1m_1)\; (j_2m_2)} $ for given two vectors 
$ \ket{j_1\; m_1}, \ket{j_2\; m_2} $. The rest steps of deriving a 
decomposition rule for \jor is the same as the case of $sl(2)$. 
We follow the standard textbook of the quantum mechanics \cite{mes}. 

  Acting $ \Delta(Z_+) $ and $ \Delta(Z_-) $ on 
$ \ket{(j_1m_1)\; (j_2m_2)} $, we obtain a series of 
vectors which are eigenvectors of $ \Delta(H) $ with eigenvalues 
\[
 -2j, \ \cdots, \ 2(m-1), \ 2m,\ 2(m+1), \ \cdots,\ 2j,
\]
where $ m=m_1+m_2 $ and $j$ denotes the highest weight. 
Let us set $N(j)$ the number of irreducible representations 
with highest weight $j$, and $n(m)$ the number of eigenvectors 
of $\Delta(H)$ with eigenvalue $2m$. The number of degenerate 
vectors can be written by the number of irreducible representations
\beq
  n(m) = \sum_{j \geq |m|}\; N(j),               \label{nirrep}
\eeq
therefore
\beq
  N(m) = n(m) - n(m+1).                          \label{Nm}
\eeq
Since $n(m)$ equals to the number of pairs $(m_1,m_2)$ 
satisfying $ m=m_1+m2 $, it can be expressed as 
\beq
 n(m) = \left\{
 \begin{array}{cl}
 0 & {\rm  for}\ \  |m| > j_1+j_2 \\
 j_1+j_2+1-|m| & {\rm for} \ \ j_1+j_2 \geq |m| \geq |j_1 - j_2| \\
 2j_2+1 & {\rm for} \ \ |j_1-j_2| \geq |m| \geq 0
 \end{array}
 \right.                    \label{states}
\eeq
Substituting (\ref{states}) into (\ref{Nm}), we obtain
\beq
 N(m) = \left\{
 \begin{array}{cl}
 1 & {\rm for} \ \ j_1+j_2 \geq |m| \geq |j_1-j_2| \\
 0 & {\rm otherwise}
 \end{array}
 \right.                    \label{numbers}
\eeq

  Therefore we have proved the fact : a tensor product of two highest 
weight representations (highest weights are $j_1$ and $j_2$) of 
\jor is reducible and the irreducible decomposition rule is 
shown schematically 
\[
 j_1 \otimes j_2 = (j_1+j_2) \oplus (j_1+j_2-1) \oplus \cdots 
 \oplus |j_1-j_2|.
\]
Furthermore each irreducible representations contained in a 
tesor product is multiplicity free.

\setcounter{equation}{0}
\section{Some examples for \jor}

  In this section, some explicit examples of irreducible decomposition, 
namely some Clebsch-Gordan coefficients, are given.
To this end, the explicit form of $\Delta(Z_-)$ is needed. Note that 
the explicit form of $\Delta(Z_+)$ is not necessary, since the vector 
which is annihilated by $\Delta(X)$ is also annihilated by $\Delta(Z_+)$. 

  From (\ref{rl}),
\beq
  \Delta(Z_-) = \Delta(\cosh \X2) \Delta(Y) \Delta(\cosh \X2). \label{copm}
\eeq
Using
\beq
 \Delta(\cosh \X2) 
  = \cosh \X2 \otimes \cosh \X2 + \sinh \X2 \otimes \sinh \X2, 
    \label{copcos}
\eeq
and (\ref{caseigen}), (\ref{ehx}), $\Delta(Z_-)$ can be rewritten as 
\bea
 & & \Delta(Z_-) \nn \\
 & &= Z_- \otimes \sum_{n=0}^{\infty}\; (n+1) \left(\Zp2 \right)^n + 
     \sum_{n=0}^{\infty}\; (n+1) \left(-\Zp2 \right)^n \otimes Z_- \nn \\
 & &+ h\left(C-{H^2 \over 4} \right) \otimes \sum_{m=1}^{\infty}\; 
     m \left(\Zp2 \right)^m - 
     \sum_{m=1}^{\infty}\; m \left( -\Zp2 \right)^m \otimes 
     h \left(C-{H^2 \over 4}\right)             \label{dzm} \\
 & &+ \left({h \over 2}\right)^2 Z_+Z_-Z_+ \otimes 
     \sum_{k=2}^{\infty} \; (k-1) \left(\Zp2 \right)^k + 
     \sum_{k=2}^{\infty} \; (k-1) \left(-\Zp2 \right)^k \otimes 
     \left({h \over 2}\right)^2 Z_+Z_-Z_+.          \nn
\eea

  We consider the cases of $m=j_1+j_2,\ j_1+j_2-1 $ and $ j_1+j_2-2$. 
Using the result of \S 3, the eigenvectors of $ \Delta(H) $ with 
eigenvalues $2m$ are constructed

\noindent
(1) $ m= j_1 + j_2 $ 
\beq
 \ket{(j_1j_1)\; (j_2j_2)} = \ket{j_1 \; j_1} \otimes 
 \ket{j_2 \; j_2},                                    \label{j1j2}
\eeq

\noindent
(2) $ m= j_1 + j_2 - 1 $
\bea
 & & \ket{(j_1j_1)\; (j_2\; j_2-1)} = \ket{j_1\; j_1} \otimes 
     \ket{j_2\; j_2-1} - hj_1 \ket{j_1; j_1} \otimes \ket{j_2 \; j_2}, 
                                                    \label{m1p1}\\
 & & \ket{(j_1\; j_1-1)\; (j_2j_2)} = \ket{j_1\; j_1-1} \otimes 
     \ket{j_2\; j_2} + 
     hj_2 \ket{j_1; j_1} \otimes \ket{j_2 \; j_2},
                                                    \label{m1p2}
\eea

\noindent
(3) $m=j_1+j_2-2$
\bea
 \ket{(j_1j_1) \; (j_2 \; j_2-2)} &=& 
 \ket{j_1\; j_1} \otimes \ket{j_2 \; j_2-2} - hj_1 
 \ket{j_1\; j_1} \otimes \ket{j_2 \; j_2-1} \nn \\
 &+& {h^2 \over 4} j_1 (2j_1-1) 
 \ket{j_1 \; j_1} \otimes \ket{j_2 \; j_2},\label{m2p1}
\eea
\bea
 & & \ket{(j_1\; j_1-1) \; (j_2 \; j_2-1)} \nn \\ 
 & & =
 \ket{j_1\; j_1-1} \otimes \ket{j_2 \; j_2-1} - h(j_1-1) 
 \ket{j_1\; j_1-1} \otimes \ket{j_2 \; j_2} \label{m2p2} \\
 & &+ h(j_2-1) \ket{j_1\; j_1} \otimes 
     \ket{j_2\; j_2-1} - {h^2 \over 2} (2j_1 j_2-j_1-j_2) 
 \ket{j_1 \; j_1} \otimes \ket{j_2 \; j_2},\nn
\eea
\bea
\ket{(j_1\; j_1-2) \; (j_2 j_2)} &=& 
 \ket{j_1\; j_1-2} \otimes \ket{j_2 \; j_2} + hj_2 
 \ket{j_1\; j_1-1} \otimes \ket{j_2 \; j_2} \nn \\
 &+& {h^2 \over 4} j_2 (2j_2-1) 
 \ket{j_1\; j_1} \otimes \ket{j_2 \; j_2},\label{m2p3}
\eea

  Let us construct the representation basis with highest weight 
$ j_1+j_2,\ j_1+j_2-1 $ and $ j_1+j2-2. $ It is easy to verify 
$ \Delta(X) \ket{(j_1j_1)\; (j_2j_2)} = 0 $ and 
$ \Delta(X) \ket{(j_1\; j_1-1)\; (j_2j_2)} = 
  \Delta(X) \ket{(j_1j_1)\; (j_2\; j_2-1)} = \ket{(j_1j_1)\; (j_2j_2)},
$ 
therefore we obtain
\bea
 & & \ket{j_1+j_2\; j_1+j_2} 
 = \ket{j_1\; j_1} \otimes \ket{j_2\; j_2},    \label{high}
\\
 & & \ket{j_1+j_2-1\; j_1+j_2-1} = \ket{(j_1\;j_1-1)\; (j_2j_2)} 
 - \ket{(j_1j_1)\; (j_2\; j_2-1)},             \label{secondhigh}
\eea
The similar calculation gives
\bea
 & & \ket{j_1+j_2-2\; j_1+j_2-2} = \nn \\
 & & \ket{(j_1j_1)\; (j_2\; j_2-2)} - \ket{(j_1\; j_1-1)\; (j_2\; j_2-1)} 
    +\ket{(j_1\; j_1-2)\; (j_2j_2)}.          \label{thirdhigh}
\eea
Other basis vectors are obtained by acting $\Delta(Z_-)$ on 
the highest weight vectors. They read
\bea
 & & \ket{j_1+j_2\; j_1+j_2-1} = {1 \over j_1+j_2} 
     (\; j_1 \ket{(j_1\; j_1-1)\; (j_2j_2)} + j_2 
          \ket{(j_1j_1)\; (j_2\; j_2-1)}\; ), \nn \\
 & & \ket{j_1+j_2\; j_1+j_2-2} = {1 \over (j_1+j_2)(2j_1+2j_2-1)} 
     \{ \; j_2(2j_2-1) \ket{(j_1j_1)\; (j_2j_2-2)} 
     \nn \\
 & & \hspace{2cm}    + 2j_1 j_2 
     \ket{(j_1\; j_1-1)\; (j_2 \; j_2-1)} + j_1(2j_1-1) 
     \ket{(j_1\; j_1-2)\; (j_2 \; j_2)} \; \},\nn \\
 & & \ket{j_1+j_2-1\; j_1+j_2-2} = {1 \over j_1+j_2-1} \{\; 
     -(2j_2-1) \ket{(j_1j_1)\; (j_2\; j_2-2)} \nn \\
 & & \hspace{2cm} 
     +(j_2-j_1) \ket{(j_1\; j_1-1)\; (j_2\; j_2-1)} +(2j_1-1) 
     \ket{(j_1\; j_1-2)\; (j_2j_2)}\; \}.   \nn 
\eea
It is remarkable that the Clebsch-Gordan coefficients for the 
vectors $ \ket{(j_1m_1)\; (j_2m_2)} $ considered in this section 
are the same as the classical ones. It may be a future work to 
investigate whether it holds for any Clebsch-Gordan coefficients.

\setcounter{equation}{0}
\section{\jorsu and its representations}

  We define \jorsu as an algebra isomorphic to \jor. Denoting the 
generators of \jorsu by $R, V, F$, they are defined
\beq
  R = -X, \qquad V = Y, \qquad F = H.           \label{defsu}
\eeq
This definition is inspired from the isomorphism between $ sl(2) $ and 
$su(1,1)$
\beq
    K_{\pm} = \mp J_{\pm}, \qquad K_0 = J_0,  \label{isom}
\eeq
where $J_{\pm}, J_0$ and $ K_{\pm}, K_0 $ are generators of 
$sl(2)$ and $su(1,1)$ respectively. Combining the isomorphism 
(\ref{isom}) and the mapping from $sl(2)$ to \jor (inverse of 
(\ref{rl})), the isomorphism (\ref{defsu}) is obtained. 

  All algebraic structure of \jorsu can easily 
be derived using (\ref{defsu}). The commutation relations are 
obtained from (\ref{com1})
\bea
 & & [F, \; R] = 2{\sinh hR \over h}, \qquad 
     [F, \; V] = -V (\cosh hR) - (\cosh hR) V, \nn \\
 & & [R, \; V] = -F,                           \label{comsu1}
\eea
the Casimir element is from (\ref{cas})
\beq
  C' = -{1 \over 2h} \{ V (\sinh hR) + (\sinh hR) V \} 
  + {1 \over 4} F^2 + {1 \over 4} (\sinh hR)^2.  \label{cassu}
\eeq
The Hopf algebra mappings for \jorsu are obtained from (\ref{hopf}).

  Let us next consider representations of \jorsu. The strategy is 
the same as the one for \jor. We define new elements of \jorsu
\beq
 T_+ = {2 \over h} \tanh \R2, \qquad 
 T_- = \left(\cosh \R2 \right) V \left( \cosh \R2 \right),
                                               \label{rlsu}
\eeq
then $ T_{\pm}, F $ satisfy the $su(1,1)$ commutation relations
\beq
  [F, \; T_{\pm}] = \pm 2 T_{\pm}, \qquad 
  [T_+,\; T_-] = - F,                          \label{comsu2}
\eeq
and the Casimir element reads
\beq
  C' = {F \over 2} \left( {F \over 2} -1 \right) - T_+ T_-. 
                                               \label{cassu2}
\eeq
These are easily verified with the identities
\bea
 & & [F,\; R^n] = 2nR^{n-1} {\sinh hR \over h}, \nn \\
 & & [V, \; R^n] = n R^{n-1} F + n(n-1) {\sinh hR \over  h}. 
                                               \label{formula}
\eea

  It is now clear that we can take the same representations 
for $ T_{\pm}, F $ as $ su(1,1) $. In this paper, we concentrate 
on the representation called the positive descrete series 
\cite{hb,ui} which is a lowest weight infinite 
dimensional representation. For the sake of simplicity of calculation, 
we adopt the different convention from \cite{hb,ui} 
\bea
 & & F \ket{\kappa \mu} = 2 \mu \ket{\kappa \mu}, \nn \\
 & & T_+ \ket{\kappa \mu} = \ket{\kappa \; \mu+1}, \label{posdes} \\
 & & T_- \ket{\kappa \mu} = (\mu - \kappa ) (\mu + \kappa - 1) 
                            \ket{\kappa \; \mu-1},  \nn
\eea
and the eigenvalu of the Casimir element is given by 
\beq
  C' \ket{\kappa \mu} = \kappa (\kappa-1) \ket{\kappa \mu}, \label{cdeigen}
\eeq
where $ \kappa $ can take any value and 
$ \mu = \kappa, \kappa +1, \kappa+2, \cdots. $ The representation 
matrices for $ R, V $ can be obtained using the inverse of 
(\ref{rlsu}).

\setcounter{equation}{0}
\section{Decomposition rule for \jorsu}

  In this section, we show that a decomposition rule of 
the product of two positive descrete series of \jorsu 
is the same as $su(1,1)$. We consider a tensor product representation 
of positive discrete series with the lowest weight $\kappa_1, 
\kappa_2. $ 
Using (\ref{rlsu}), the coproduct of $F$ can be rewritten as 
\beq
  \Delta(F) = F \otimes 1 + 1 \otimes F   \nn \\
            + F \otimes 2 \sum_{n=1}^{\infty}\; 
                  \left(-\Tp2 \right)^n + 
                  2 \sum_{n=1}^{\infty}\; 
                  \left(\Tp2 \right)^n \otimes F. \label{copf}
\eeq
The eigenvector of $\Delta(F)$ with eigenvalu $ 2(\mu_1 + \mu_2) $ 
may be written
\beq
  \ket{(\kappa_1\mu_1)\; (\kappa_2\mu_2)} = 
  \sum_{\rho, \sigma = 0}^{\infty}\; 
  \alpha(\mu_1+\rho,\mu_2+\sigma) \; 
  \ket{\kappa_1\; \mu_1+\rho} \otimes 
  \ket{\kappa_2\; \mu_2+\sigma}.          \label{eigensu}
\eeq
Because of the consistency with the limit of $ h \longrightarrow 0 $, 
we set $ \alpha(\mu_1, \mu_2) = 1. $ Substituting (\ref{copf}) 
and (\ref{eigensu}) into the relation 
$ \Delta(F) \ket{(\kappa_1\mu_1)\; (\kappa_2\mu_2)} = 
2 (\mu_1 + \mu_2) \ket{(\kappa_1\mu_1)\; (\kappa_2\mu_2)} $, 
we obtain the recurrence relation for 
$ \alpha(\mu_1+\rho, \mu_2+\sigma) $
\bea
 & & (\rho + \sigma) \alpha(\mu_1+\rho, \mu_2+\sigma) 
     + 2 (\mu_1 + \rho) \sum_{n=1}^{\sigma}\; \left(-{h \over 2}\right)^n \; 
       \alpha(\mu_1+\rho, \mu_2+\sigma-n)      \nn \\
 & & + 2 (\mu_2 + \sigma) \sum_{n=1}^{\rho}\; \left({h \over 2} \right)^n \; 
       \alpha(\mu_1+\rho-n, \mu_2+\sigma) = 0.  \label{recsu1}
\eea
Repeating the same procedure as the case of \jor, the recurrence 
relation (\ref{recsu1}) is rewritten into the simpler form
\bea
 & & (\rho + \sigma) \; \alpha(\mu_1+\rho, \mu_2+\sigma) 
     - {h \over 2} (2\mu_1+1+\rho-\sigma) 
     \; \alpha(\mu_1+\rho, \mu_2+\sigma -1) \nn \\
 & & + {h \over 2} (2 \mu_2+1-\rho+\sigma) 
     \; \alpha(\mu_1+\rho-1, \mu_2+\sigma)  \label{recsu2} \\
 & &  + 
     \left( {h \over 2} \right)^2 (2\mu_1+2\mu_2-2+\rho+\sigma) 
     \; \alpha(\mu_1+\rho-1, \mu_2+\sigma -1) = 0.   \nn
\eea

  The solutions of (\ref{recsu2}) are given by 
\bea
 & & \alpha(\mu_1+\rho, \mu_2 + \sigma) \nn \\
 & & = (-1)^{\rho} \left( {h \over 2} \right)^{\rho + \sigma} 
 \sum_{p=0}\; 
 \left( \begin{array}{c}
        2\mu_1 + \rho - p \\ \sigma - p
        \end{array} \right) 
 \left( \begin{array}{c}
        2\mu_1 + \rho - 1 \\ p
        \end{array} \right) 
 \left( \begin{array}{c}
        2\mu_2 \\ \rho - p
        \end{array} \right)_, \label{solsu}
\eea
where the sum on $p$ runs as far as all the binomial coefficients 
are welldefined. It can be proved that (\ref{solsu}) satisfies 
the recurrence relation (\ref{recsu2}) in the same way as 
\S 3.

  It has been shown that we can construct unique eigenvector of 
$ \Delta(F) $ with eigenvalu $ 2(\mu_1+\mu_2) $ for given 
vectors $ \ket{\kappa_1 \mu_1}, \ket{\kappa_2 \mu_2}. $ 
Acting $ \Delta(T_{\pm}) $ on 
$ \ket{(\kappa_1 \mu_1)\; (\kappa_2\; \mu_2)}, $ 
we can construct a series of eigenvectors of $ \Delta(F) $ with 
eigenvalues 
\[
 \kappa, \ \kappa + 1, \ \cdots, \ \mu, \ \mu +1, \ \cdots,
\]
where $ \mu = \mu_1 + \mu_2 $ and it is clear that the 
lowest possible value of $ \mu $ (denoted by $ \kappa $) is 
$ \kappa_1 + \kappa_2. $ Let us set $ N(\kappa ) $ the number 
of irreducible representations with lowest weight $ \kappa $, 
and $ n(\mu) $ the number of eigenvectors of $ \Delta(F) $ with 
eigenvalu $ 2\mu $. The number of degenerate vectors can be 
written by the number of irreducible representation
\beq
  n(\mu) = \sum_{\kappa \leq \mu}\; N(\kappa),   \label{nmu}
\eeq
therefore
\beq
  N(\mu) = n(\mu) - n(\mu -1).                  \label{lnmu}
\eeq
Since $ n(\mu) $ equals the number of pairs $ (\mu_1, \mu_2) $ 
satisfying $ \mu = \mu_1 + \mu_2 $, it is given by 
\beq
  n(\mu) = \left\{
  \begin{array}{cl}
  0 & {\rm for} \quad \mu < \kappa_1 + \kappa_2 \\
  \mu - \kappa_1 - \kappa_2 + 1 & 
  {\rm for} \quad \mu \geq \kappa_1+ \kappa_2
  \end{array}
  \right.                                      \label{pair}
\eeq
Substituting (\ref{pair}) into (\ref{lnmu}), 
\beq
  N(\mu) = \left\{
  \begin{array}{ll}
  0 & {\rm for} \quad \mu < \kappa_1 + \kappa_2 \\
  1 & 
  {\rm for} \quad \mu \geq \kappa_1+ \kappa_2
  \end{array}
  \right.                                      \label{lnmu2}
\eeq

 Therefore we have proved the fact : a tensor product of two 
positive descrete series of \jorsu is reducible and the irreducible 
decomposition rule is given schematically by 
\[
  \kappa_1 \otimes \kappa_2 = \kappa_1 + \kappa_2, \ 
  \kappa_1 + \kappa_2+1, \ \kappa_1 + \kappa_2 + 2, \ \cdots.
\]
Furthermore each irreducible representation contained in the 
tensor product is multiplicity free.

\setcounter{equation}{0}
\section{Conclusion}

  We have shown that, for both highest weight finite dimensional 
representations of \jor and lowest weight infinite dimensional ones 
of \jorsu, tensor product representations are reducible and 
the decomposition rules to irreducible representations are 
exactly the same as those of the corresponding Lie algebras. 
We concentrate on the positive descrete series of \jorsu, the 
same result may hold for the negative descrete series which 
are highest weight infinite dimensional representation, since 
the difference between positive and negative descrete series is 
to use highest weight or lowest one. The Lie algebra $ su(1,1) $ has 
two other infinite dimensional representations \cite{hb}. 
The corresponding representations of \jorsu may obtain the 
inverse mapping of (\ref{rlsu}), however tensor products of 
such representations are still an open problem. 

  The construction of eigenvectors of $ \Delta(H) $ and $ \Delta(F) $ 
is the key of the proof. The other steps of the proof are nothing 
but the ones for the Lie algebras. These parallelism in the 
representation theories between Jordanian quantum algebras 
and the corresponding Lie algebras may suggest further similarities.
For example, we might be able to obtain the Clebsch-Gordan 
coefficients by the same method as the classical case, Racha-Wigner 
type of calculus (6-j, 9-j symbols, tensor operators, 
Wigner-Eckart's theorem {\em etc}) might be possible for the Jordanian 
quantum algebras. The similarity in the representation theories may 
also suggest that the Jordanian quantum algebra are applicable to 
various fields in physics. These will be future works.

%

\end{document}